\documentclass{article}
\usepackage{pst-node}
\usepackage{arxiv}
\usepackage[english]{babel}
\usepackage[utf8]{inputenc} 
\usepackage[T1]{fontenc}    
\usepackage{hyperref}
\hypersetup{
    colorlinks=true,
    linkcolor=red,
    urlcolor=blue,
    citecolor=black}
\usepackage{url}            
\usepackage{booktabs}       
\usepackage{amsfonts}       
\usepackage{nicefrac}       
\usepackage{microtype}      
\usepackage{lipsum}
\usepackage{amsmath}
\usepackage{amssymb}
\usepackage{pdfpages}
\usepackage[square,numbers]{natbib}
\usepackage[nottoc]{tocbibind}
\usepackage{natbib}
\usepackage{tikz}
\usepackage{graphicx}
\usetikzlibrary{positioning}

\usepackage{url}            
\usepackage{amsfonts}       
\usepackage{graphicx}

\title{Imaginary Cities of the Diophantine equation $X^3+Y^3+Z^3=K$}

\author{
  Eduardo J. Acu{\~n}a T.\\
  \texttt{eacuna1@uc.edu.ve} \\
  Universidad de Carabobo
 \\
 Facultad de Ciencia y Tecnologia.
}

\begin{document}
\maketitle

\begin{abstract}
    The purpose of this article is to make a graph representation of the Diophantine equation $X^3+Y^3+Z^3=K$ using the theory of "imaginary cities" for its construction and to determine the modular combinations on its edges with the De Bruijn cycles.
\end{abstract}

\section{Introduction}
the famous Diophantine equation $X^3+Y^3+Z^=K$ better known as the sum of the three cubes, has long proved to be a challenging problem, especially with the two questions posed by this problem\\
How many integer $K$ solutions does this equation have, or put another way how many ways can $K$ be written as the sum of 3 cubes?\\
We can determine that $K$ has no solution for numbers that when dividing $K/9$ the remainder is 4 or 5, i.e. it is not a solution for all integers $K = +-4$ mod 9.\\
Heath-Brown \cite{heath1992density} makes the conjecture that there are infinitely many solutions to $K$ as long as it satisfies the above condition for $K$ to be a solution.\\
The problem is that for some small solutions, the cubes can be immensely large and that is what makes this problem difficult, however in 2019 Andrew Booker \cite{booker2019cracking} manages to find the solution to the number 33, subsequently, in 2021\cite{booker2021question} publishes the solution to 42 using a computational algorithm, as well as the algorithm left by Elkies \cite{elkies2000rational}(2000) and later Huisman \cite{huisman2016newer}(2016) who left us extensive data.\\
In 2021 \cite{romero2021matrixe_phi} a matrix-based algebra is proposed that provides a modular class system and studies its behavior between operations, realizing that mods can be the solution for this type of problem it was decided to use this theory in the three-cube problem, later that year a pre-print \cite{flores2021existential} is published with a theoretical and non-computational algorithm, showing that for each solution $K=z$ mod 9, depending on the z value has finite paths to find those solutions.\\
With graph theory and specifically using De Bruijn Cycles we can determine that the same finite paths are within a numerical cycle forming a graph that represents all possible combinations and therefore all possible solutions that exist for this problem.

\section{The mod is the key}
The first clue we have is the mod, to know which numbers are not the result of this equation, therefore we must look for other numbers within the remaining mod and see if a number is not a result so we would say that the equation does not have infinite $K$ solutions in those mods.\\

Using a modular analysis as done in \cite{romero2021matrixe_phi} we determine that there exists a set of classes for each integer in mod 9, where the remainder $z/9$ will be its class, in that article, for definition purposes, the numbers congruent to 0 mod 9 were called class 9, since all of them are multiples of 9, for this case, as it is desired to simplify the class, it will take the value only of its remainder in that mod.

\begin{equation*}
    n\in \mathbb{Z}~|~n/9~mod~9\cong z=\{0,..,8\}
\end{equation*}

It is already trivially known that cubed integers always have the following congruence's in mod 9:

\begin{equation*}
    \frac{(n)^3}{9}\cong\{0,1,8\}~(mod~9)
\end{equation*}

This means that we can rewrite the equation in its modular form, and the result that it gives us is the modulus of the solution of $K$.\\
It is evident from what was explained in this pre-print \cite{flores2021existential} that the mods of $\phi^K$ with the most difficult solutions to find are 3 and 6.

 \subsection*{case 1:}With everything explained, we can carry out the following case, if we wanted to find the modular path that gives us all the $\phi^K=6$ mod 9 solutions, they would be the following:\\

\begin{center}   
$\begin{array}{cc}
\label{a1}
     8+8+8&=6  \\
     8-1-1&=6 \\
     8+8-1&=6
\end{array}$
\end{center}

This same combination can be done to obtain all possible mod that are solutions to this equation.

\subsection*{example:}Let's do an example with the combination 8-1-1=6 for that I will take numbers from the Huisman data in his article\cite{huisman2016newer} and I will verify that they all meet the unique combinations shown.\\
\begin{itemize}
    \item[step 1.] We determine that the number we want to find has the same congruence with the path that we are going to implement.\\
      $K=15  \cong \phi^k=6$
    \item[step 2.] We convert the numbers to their modular form without raising them to the cube.\\
    $K=(-265)^3 + (-262)^3 + (332)^3  \cong (-4_\phi)^3+(-1_\phi)+(8_\phi)$
    \item[step 3.]We cube all the modular forms and we will also take the class according to the rest of the results.\\
     $(-4_\phi)^3+(-1_\phi)+(8_\phi) \cong (-1)+(-1)+8$
    \item[step 4.] We verify that the combination obtained is congruent with the result.\\
       $8+(-1)+(-1)=6 \cong (-265)^3 + (-262)^3 + (332)^3=15$

\end{itemize}
The famous result 33 was also found in this combination.

\section{From mod to graph}
First let's start by defining a graph as a set of nodes and edges that can represent maps, knots, paths, and many other things, in this case, we will use it to represent a sum of three cubes.
\begin{equation*}
    G=(n,e)
\end{equation*}
where, $n$ is a node, and $e$ is the edge it is attached to.\\
There are a wide variety of ways to define a graph type, but we want to determine if this equation is infinitely cyclic, so it should be able to be described in an Eulerian graph.\\

To find the graph form of $X^3+Y^3+Z^3$ we will use the "imaginary city" I.J. Good.\cite{good1946normal} (1946) to generate a De Bruijn sequence, as we can see in figure 5, $N_2$ in the \emph{"A combinatorial problem. Proceedings of the Section of Sciences of the Koninklijke Nederlandse Akademie van Wetenschappen te Amsterdam, 49(7). 758-764.}\cite{de1946combinatorial}:

\subsection*{example 2:}
let's make a trip from island 00 to island 01, to name the bridge that leads us we will take 00 and intercept it with 01, and at the end, we will get 001 which would be our edge:

\begin{equation*}
    N(00)\rightarrow N(01),~ \underbrace{00}_{begin}\cap\overbrace{01}^{end}=(001)
\end{equation*}
In the case of those that intercept themselves, we have $E_0(000)$ and $E_1(111)$. If we take all the values of the edges and intercept them side by side we will find a sequence and these are what we call "De Bruijn cycles", the ones that generated the graph. The one that generated the graph is this one: $P_3~cycle=(00010111)$ which contains all the triples that give names to the edges, (000),(001),(010),(101),(011),(111).

\section{\texorpdfstring{Graph of $X^3+Y^3+Z^3$}{}}

\subsection{Graph structure:}
\begin{enumerate}
    \item Nodes: The nodes are composed of the first two or the last two of the triads.
    \item Edges: show the unique address in which to combine the triads and will be named after the combination of the nodes they link
    \item Dashed edges: are the edges that have already been seen show where the cycle begins to repeat itself infinitely.
\end{enumerate}
\subsection{Graphs}
Using the same principles that we have just seen, this time we will build an imaginary city with edges that are also triples but instead of using a binary system, we will use a ternary system, since the digits that we are going to use are the only three digits that are obtained as remainder in mod 9 of a cube, (0,1,8).\\
When a ternary system is used to generate the bridges, we realized that what is happening is that "imaginary cities" are being formed that are interconnecting with each other continuously and infinitely.

\begin{center}
    \begin{tikzpicture}
[g/.style={circle, draw=black!60, fill=white!5, very thick, minimum size=7mm},
a/.style={circle, draw=black!60, fill=white!5, very thick, minimum size=7mm},
r/.style={circle, draw=black!60, fill=white!5, very thick, minimum size=7mm},
n/.style={circle, draw=white!60, fill=white!5, very thick, minimum size=7mm},
]

\node[n]               (n1)                     {};
\node[a]               (11)    [left=of n1]       {11};
\node[a]               (00)    [right=of n1]       {00};
\node[g]               (10)    [above=of n1]       {10};
\node[r]               (01)    [below=of n1]       {01};
\node[n]               (n2)    [right=of 00]       {};
\node[g]               (08)    [above=of n2]       {08};
\node[r]               (80)    [below=of n2]       {80};
\node[a]               (88)    [right=of n2]       {88};
\node[n]               (n3)    [right=of 88]       {};
\node[g]               (81)    [above=of n3]       {81};
\node[r]               (18)    [below=of n3]       {18};
\node[a]               (111)    [right=of n3]       {11};
\node[n]               (n4)    [right=of 111]       {};
\node[g]               (110)    [above=of n4]       {10};
\node[r]               (011)    [below=of n4]       {01};
\node[n]               (n5)    [right=of n4]       {};

\draw[->, ultra thick] (11.north) -- (10.west);
\draw[->, ultra thick] (10.east)  -- (08.west);
\draw[->, ultra thick] (08.east) -- (81.west);
\draw[dashed, ->, ultra thick] (81.east) -- (110.west);
\draw[->, ultra thick] (01.north) .. controls +(down:1mm) and +(right:9mm) .. (10.south);
\draw[->, ultra thick] (10.south) .. controls +(down:1mm) and +(left:9mm) .. (01.north);
\draw[->, ultra thick] (80.north) .. controls +(down:1mm) and +(right:9mm) .. (08.south);
\draw[->, ultra thick] (08.south) .. controls +(down:1mm) and +(left:9mm) .. (80.north);
\draw[->, ultra thick] (18.north) .. controls +(down:1mm) and +(right:9mm) .. (81.south);
\draw[->, ultra thick] (81.south) .. controls +(down:1mm) and +(left:9mm) .. (18.north);
\draw[dashed, ->, ultra thick] (011.north) .. controls +(down:1mm) and +(right:9mm) .. (110.south);
\draw[dashed, ->, ultra thick] (110.south) .. controls +(down:1mm) and +(left:9mm) .. (011.north);
\draw[->, ultra thick] (01.west) -- (11.south);
\draw[->, ultra thick] (80.west)  -- (01.east);
\draw[->, ultra thick] (18.west) -- (80.east);
\draw[->, ultra thick] (011.west) -- (18.east);
\draw[->, ultra thick] (10.east) -- (00.north);
\draw[->, ultra thick] (00.south)  -- (01.east);
\draw[->, ultra thick] (08.east) -- (81.west);
\draw[->, ultra thick] (08.east) -- (88.north);
\draw[->, ultra thick] (88.south)  -- (80.east);
\draw[->, ultra thick] (80.west)  -- (00.south);
\draw[->, ultra thick] (00.north)  -- (08.west);
\draw[dashed, ->, ultra thick] (81.east)  -- (111.north);
\draw[dashed, ->, ultra thick] (111.south)  -- (18.east);
\draw[->, ultra thick] (18.west)  -- (88.south);
\draw[->, ultra thick] (88.north)  -- (81.west);
\draw[dashed, ->, ultra thick] (011.west)  -- (111.south);
\draw[dashed, ->, ultra thick] (111.north)  -- (110.west);
\draw[dashed, ->, ultra thick] (110.east)  -- (n5.north);
\draw[dashed, ->, ultra thick] (n5.south)  -- (011.east);

\end{tikzpicture}

\label{Fig1}
\textbf{Figure 1}

\end{center}
Although this continuous graph shows me that the sequence formed by the three mod in all their possible combinations extends infinitely, it does not show me if it is cyclical, although it seems evident.\\
In the following graph, we manage to unite all the nodes with their correspondents in this way we obtain all the edge triplets in a single drawing.

\begin{center}
    
\begin{tikzpicture}[
G8/.style={circle, draw=black!60, fill=white!5, very thick, minimum size=7mm},
G0/.style={circle, draw=black!60, fill=white!5, very thick, minimum size=7mm},
G1/.style={circle, draw=black!60, fill=white!5, very thick, minimum size=7mm},
In/.style={circle, draw=white!60, fill=white!5, very thick, minimum size=7mm},
]
\node[In]               (in1)                     {};
\node[G0]               (88)       [left=2cm of in1]  {88};
\node[G0]               (00)       [right=2cm of in1] {00};
\node[G8]               (08)      [below=of in1]   {08};
\node[In]               (in2)       [left=of 08] {};
\node[In]               (in3)       [below=of in2] {};
\node[G1]               (18)       [left=of in3] {18};
\node[In]               (in4)       [right=of 08] {};
\node[In]               (in5)       [below=of in4] {};
\node[G1]               (01)       [right=of in5] {01};
\node[In]               (in6)       [below=of 18] {};
\node[In]               (in7)       [below=of 01] {};
\node[G0]               (11)       [left=2cm of in6] {11};
\node[G0]               (111)       [right=2cm of in7] {11};
\node[G8]               (10)       [below=of in6] {10};
\node[G8]               (81)       [below=of in7] {81};
\node[In]               (in8)       [right=of 10] {};
\node[In]               (in9)       [below=of in8] {};
\node[G1]               (80)       [right=of in9] {80};
\node[In]               (invi)       [below=of 80] {};
\node[G0]                (000)      [left=2cm of invi] {00};
\node[G0]               (888)       [right=2cm of invi] {88};

\draw[->, ultra thick] (08.south) to node[right] {$081$} (81.north);
\draw[->, ultra thick] (81.west) to node[below] {$810$} (10.east);
\draw[->, ultra thick] (10.north) to node[left] {$108$} (08.south);
\draw[->, ultra thick] (18.south) to node[right] {$180$} (80.north);
\draw[->, ultra thick] (80.north) to node[left] {$801$} (01.south);
\draw[->, ultra thick] (01.west) to node[below] {$018$} (18.east);
\draw[->, ultra thick] (08.north) to node[above] {$088$} (88.east);
\draw[->, ultra thick] (18.north) to node[left] {$188$} (88.south);
\draw[->, ultra thick] (00.west) to node[above] {$008$} (08.north);
\draw[->, ultra thick] (00.south) to node[right] {$001$} (01.north);
\draw[->, ultra thick] (11.north) to node[above] {$118$} (18.west);
\draw[->, ultra thick] (11.south) to node[below] {$110$} (10.west);
\draw[->, ultra thick] (01.east) to node[above] {$011$} (111.north);
\draw[->, ultra thick] (81.east) to node[below] {$811$} (111.south);
\draw[->, ultra thick] (888.north) to node[right] {$881$} (81.south);
\draw[->, ultra thick] (888.west) to node[below] {$880$} (80.south);
\draw[->, ultra thick] (80.south) to node[below] {$800$} (000.east);
\draw[->, ultra thick] (10.south) to node[right] {$100$} (000.north);
\draw[dashed,<->, ultra thin] (000.north) -- (00.south);
\draw[dashed,<->, ultra thin] (888.north) -- (88.south);
\draw[dashed,<->, ultra thin] (111.west) -- (11.east);
\draw[->] (888.south) to node[below]{$888$} (888.south);
\draw[->] (000.south) to node[below]{$000$} (000.south);
\draw[->] (111.east) to node[right]{$111$} (111.east);
\draw[->] (11.west) to node[left]{$111$} (11.west);
\draw[->] (88.north) to node[above]{$888$} (88.north);
\draw[->] (00.north) to node[above]{$000$} (00.north);
\end{tikzpicture}
\label{fig2}

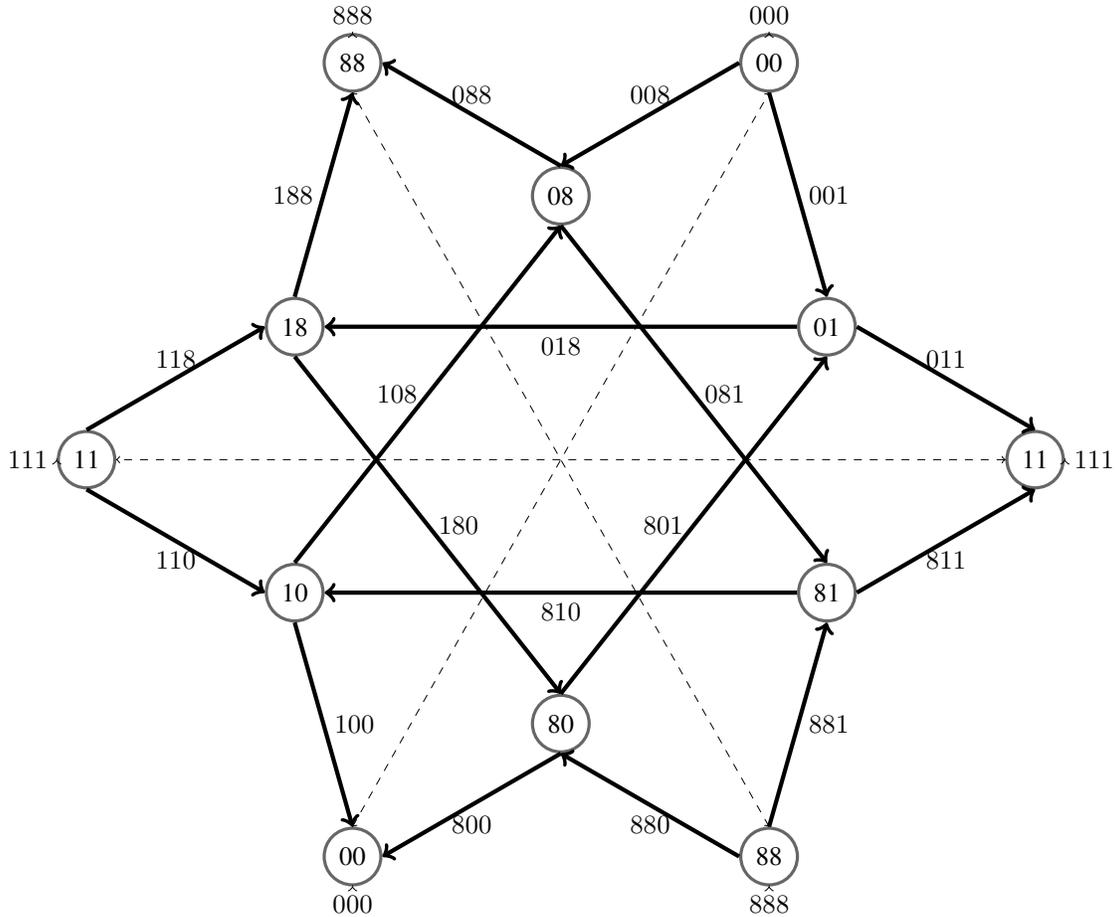
\textbf{Figure 2: Graph $X^3+Y^3+Z^3$}
\end{center}
Now that we have this graph, we can perform the same procedure and overlap all the edges to obtain the De Bruijn cycle:
\begin{equation*}
    A_3-Cycle=(00088808881118100010110)
\end{equation*}
As we know, all the edges are contained in this sequence, but the most interesting thing is that all the combinations of sums of three cubes in the modular form are also here.
\subsubsection*{Example 3:}
We can search in which part of the sequence the famous number 33 of Andrew Booker\cite{booker2019cracking} was found.
\begin{equation*}
    A_3-Cycle=(000888088\underbrace{811}_{(8+(-1)+(-1)=6}18100010110)
\end{equation*}
But it is not exclusive to that area, the result can also be found in combinations with the same digits in the cycle, showing that there are more combinations for the same result, which can be seen in the edges, but to simplify if we join all the edges that share the same elements we will have that for that specific result we have a specific path.\\
This is closely related to the Primordial Conjecture Acu{\~n}a. E. Marrero P. (2022)\cite{https://doi.org/10.5281/zenodo.6595095} that he makes: if the operation of the remainders of the integers in a mod 9, results in a congruent remainder, for said equation there are solutions for that equation with that remainder.
\newpage
\subsection{Sub-Graphs}
In the center of \textbf{Figure 2} there are also node connections, which we can visualize well in \textbf{Figure 1}, for lack of space and to make it look more orderly, it is shown below as a sub-graph, the edges that are given in this sub-graph are complementary to the sequence (The differentiation of colors in the nodes and edges, is only to improve the visualization).\\
\begin{center}
    
\begin{tikzpicture}[
G8/.style={circle, draw=black!60, fill=white!5, very thick, minimum size=7mm},
G0/.style={circle, draw=blue!60, fill=white!5, very thick, minimum size=7mm},
G1/.style={circle, draw=red!60, fill=white!5, very thick, minimum size=7mm},
In/.style={circle, draw=white!60, fill=white!5, very thick, minimum size=7mm},
]
\node[In]               (in1)                     {};

\node[G8]               (08)      [below=of in1]   {08};
\node[In]               (in2)       [left=of 08] {};
\node[In]               (in3)       [below=of in2] {};
\node[G0]               (18)       [left=of in3] {18};
\node[In]               (in4)       [right=of 08] {};
\node[In]               (in5)       [below=of in4] {};
\node[G1]               (01)       [right=of in5] {01};
\node[In]               (in6)       [below=of 18] {};
\node[In]               (in7)       [below=of 01] {};
\node[G1]               (10)       [below=of in6] {10};
\node[G0]               (81)       [below=of in7] {81};
\node[In]               (in8)       [right=of 10] {};
\node[In]               (in9)       [below=of in8] {};
\node[G8]               (80)       [right=of in9] {80};
\node[In]               (invi)       [below=of 80] {};

\draw[dashed, ->] (08.south) -- (81.north);
\draw[dashed, ->] (81.west) -- (10.east);
\draw[dashed, ->] (10.north) -- (08.south);
\draw[dashed, ->] (18.south) -- (80.north);
\draw[dashed, ->] (80.north) -- (01.south);
\draw[dashed, ->] (01.west) -- (18.east);
\draw[color=blue!60,->] (18.south) to node[below] {$181$} (81.south);
\draw[color=blue!60,->] (81.north) to node[above] {$818$} (18.north);
\draw[color=red!60,ultra thin,->] (10.south) to node[above] {$101$} (01.south);
\draw[color=red!60,->] (01.north) to node[below] {$010$} (10.north);
\draw[->] (08.west) to node[left] {$080$} (80.west);
\draw[->] (80.east) to node[right] {$808$} (08.east);
\end{tikzpicture}

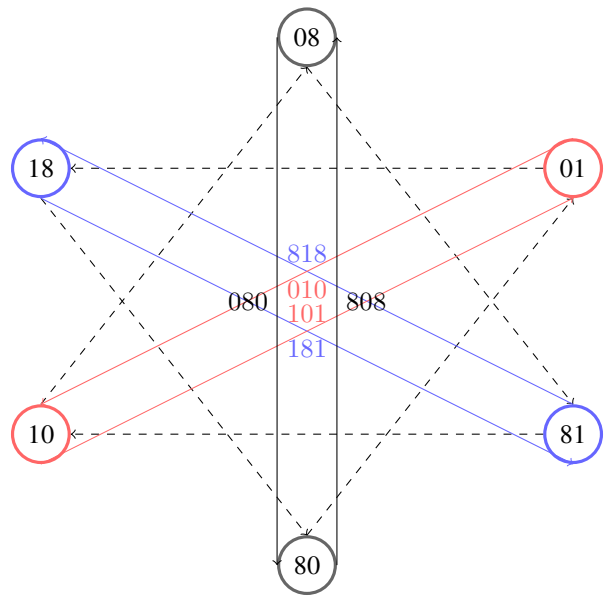
\textbf{Figure 3: Sub-graph $G_0=(E_0,N_0)$}
\end{center}

Once we have this graph \textbf{Figure 2.} we can decompose it to obtain two sub-graphs that give us more information since both sub-graphs are Eurelian\cite{bollobas1979graphs}:
\begin{center}

\begin{tikzpicture}[
g/.style={circle, draw=black!60, fill=white!5, very thick, minimum size=7mm},
a/.style={circle, draw=black!60, fill=white!5, very thick, minimum size=7mm},
r/.style={circle, draw=black!60, fill=white!5, very thick, minimum size=7mm},
n/.style={circle, draw=white!60, fill=white!5, very thick, minimum size=7mm},
]

\node[n] (n1)                  {};
\node[a] (11)  [right=0.5cm of n1] {11};
\node[n] (n2) [right=0.5cm of 11]    {};
\node[g] (18) [below=of n1]  {18};
\node[g] (01) [below=of n2]  {01};
\node[n] (n3) [below=of 18]    {};
\node[n] (n4) [below=of 01]    {};
\node[a] (88) [left=of n3]   {88};
\node[g] (80) [below=2.75cm of 11]  {80};
\node[a] (00) [right=of n4]  {00};

\draw[->, ultra thick] (11.south) to node[left] {$118$}  (18.north);
\draw[->, ultra thick] (01.north) to node[right] {$011$}  (11.south);
\draw[->, ultra thick] (01.west) to node[above] {$018$}  (18.east);
\draw[->, ultra thick] (18.west) to node[left] {$188$}  (88.north);
\draw[->, ultra thick] (00.north) to node[right] {$001$}  (01.east);
\draw[->, ultra thick] (88.east) to node[below] {$880$}  (80.west);
\draw[->, ultra thick] (80.east) to node[below] {$800$} (00.west);
\draw[->, ultra thick] (18.south) to node[left] {$180$} (80.north);
\draw[->, ultra thick] (80.north) to node[right] {$801$} (01.south);
\draw[->,thick]  (11.north) arc (1:270:4mm) node[left]{$111$};
\draw[->,thick]  (88.south) arc (360:90:4mm) node[left]{$888$};
\draw[->,thick]  (00.west) arc (90:360:4mm) node[right]{$000$};
\end{tikzpicture}

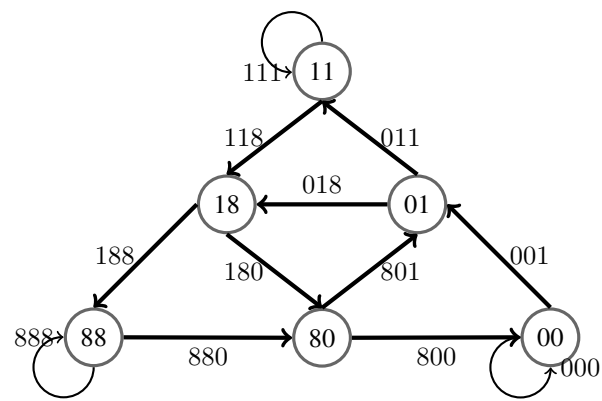
\textbf{Figure 4: Sub-graph $G_1=(E_1,N_1)$}
\end{center}
$E_1=$(011)(111)(118)(188)(888)(880)(800)(000)(001)(018)(180)(801)\\
$G_1-Cycle=$0111818880800018180801

\begin{center}

\begin{tikzpicture}[
g/.style={circle, draw=black!60, fill=white!5, very thick, minimum size=7mm},
a/.style={circle, draw=black!60, fill=white!5, very thick, minimum size=7mm},
r/.style={circle, draw=black!60, fill=white!5, very thick, minimum size=7mm},
n/.style={circle, draw=white!60, fill=white!5, very thick, minimum size=7mm},
]

\node[n] (n1)                  {};
\node[a] (11)  [right=0.5cm of n1] {11};
\node[n] (n2) [right=0.5cm of 11]    {};
\node[g] (81) [below=of n1]  {81};
\node[g] (10) [below=of n2]  {10};
\node[n] (n3) [below=of 18]    {};
\node[n] (n4) [below=of 01]    {};
\node[a] (88) [left=of n3]   {88};
\node[g] (08) [below=2.75cm of 11]  {08};
\node[a] (000) [right=of n4]  {00};

\draw[->, ultra thick] (11.south) to node[right] {$110$}  (10.north);
\draw[->, ultra thick] (81.north) to node[left] {$811$}  (11.south);
\draw[->, ultra thick] (81.east) to node[above] {$810$}  (10.west);
\draw[->, ultra thick] (88.north) to node[left] {$881$}  (81.west);
\draw[->, ultra thick] (10.east) to node[right] {$100$}  (00.north);
\draw[->, ultra thick] (08.west) to node[below] {$088$}  (88.east);
\draw[->, ultra thick] (00.west) to node[below] {$008$} (08.east);
\draw[->, ultra thick] (08.north) to node[left] {$081$} (81.south);
\draw[->,thick]  (11.north) arc (1:270:4mm) node[left]{$111$};
\draw[->,thick]  (88.south) arc (360:90:4mm) node[left]{$888$};
\draw[->, ultra thick] (10.south) to node[right] {$108$} (08.north);
\draw[->,thick]  (00.west) arc (90:360:4mm) node[right]{$000$};
\end{tikzpicture}

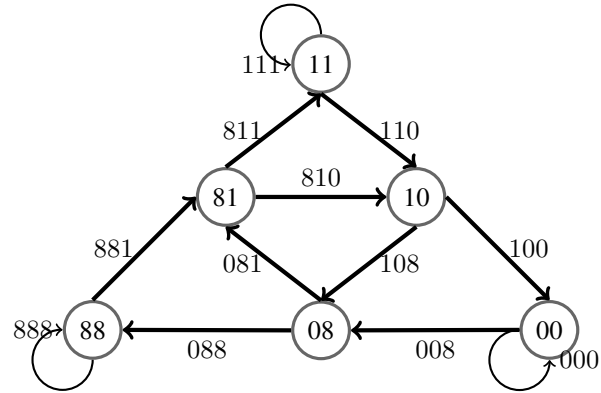
\textbf{Figure 5: Sub-graph $G_2=(E_2,N_2)$}
\end{center}
$E_2=$(811)(111)(110)(100)(000)(008)(088)(888)(881)(810)(108)(081)\\
$G_2-Cycle=$8111010008088818101081\\

As we can see each of these sub-graphs forms Eurelian cycles and De Bruijn cycles, what we are representing is that the equation goes continuously and infinitely in all its possible combinations, in terms of Euler we can cross all the bridges of them without repeating once, and this process can be repeated infinitely, with this we come to the same conclusion of Heat Brown, this equation has solutions in all its possible $K$ and also has infinite ways to be solved.

\section{Conclusion}
We can conclude that there is a strong relationship between the modular forms of an equation and the graphs.\\
In a problem of such high difficulty as this one, it is necessary to extend more and more the methods to analyze it, it has been many years and many mathematicians have contributed to being able to arrive at its solution, what we propose is a perspective from the constructionist point of view, where joining each theorem that has been given in the past and showing this visual representation we can conclude that this equation is in an infinite cycle of solutions, unlike the mythical Fermat's last theorem, this equation allows the use of negative numbers which broadens its permissiveness of solutions, therefore we conclude that it does have solutions in all its possible $K$.\\
We hope that in the future a more rigid proof will be given and that this research will make a computational and mathematical contribution to reach its final solution.

\section*{Acknowledgments}
During the years of research I have received support and comments from many mathematicians, first of all to Professor Samuel Flores, Andres Lugo, mathematician Paul Marrero,for taking some of their busy time and answering my questions, thank you very much.
\bibliography{sample}

\begin{thebibliography}{11}
\providecommand{\natexlab}[1]{#1}
\providecommand{\url}[1]{\texttt{#1}}
\expandafter\ifx\csname urlstyle\endcsname\relax
  \providecommand{\doi}[1]{doi: #1}\else
  \providecommand{\doi}{doi: \begingroup \urlstyle{rm}\Url}\fi

\bibitem[Bollob{\'a}s and Bollob{\'a}s(1979)]{bollobas1979graphs}
B.~Bollob{\'a}s and B.~Bollob{\'a}s.
\newblock Graphs and groups.
\newblock \emph{Graph Theory: An Introductory Course}, pages 146--174, 1979.

\bibitem[Booker(2019)]{booker2019cracking}
A.~R. Booker.
\newblock Cracking the problem with 33.
\newblock \emph{Research in Number Theory}, 5:\penalty0 1--6, 2019.

\bibitem[Booker and Sutherland(2021)]{booker2021question}
A.~R. Booker and A.~V. Sutherland.
\newblock On a question of mordell.
\newblock \emph{Proceedings of the National Academy of Sciences}, 118\penalty0 (11):\penalty0 e2022377118, 2021.

\bibitem[De~Bruijn(1946)]{de1946combinatorial}
N.~G. De~Bruijn.
\newblock A combinatorial problem.
\newblock \emph{Proceedings of the Section of Sciences of the Koninklijke Nederlandse Akademie van Wetenschappen te Amsterdam}, 49\penalty0 (7):\penalty0 758--764, 1946.

\bibitem[Eduardo J.~Acu{\~n}a and Paul F.~Marrero(2021)]{https://doi.org/10.5281/zenodo.6595095}
T.~Eduardo J.~Acu{\~n}a and R.~Paul F.~Marrero.
\newblock Primal conjecture in matrixe9$\phi$, 2021.
\newblock URL \url{https://zenodo.org/record/6595095}.

\bibitem[Elkies(2000)]{elkies2000rational}
N.~D. Elkies.
\newblock Rational points near curves and small nonzero| x 3- y 2| via lattice reduction.
\newblock In \emph{Algorithmic Number Theory: 4th International Symposium, ANTS-IV Leiden, The Netherlands, July 2-7, 2000. Proceedings 4}, pages 33--63. Springer, 2000.

\bibitem[Flores et~al.(2021)Flores, Acu{\~n}a, and Marrero]{flores2021existential}
S.~Flores, E.~Acu{\~n}a, and P.~Marrero.
\newblock Existential refinement on the search of integer solutions for the diophantine equation $ x^3+y^3+z^3= n$.
\newblock \emph{arXiv preprint arXiv:2103.17037}, 2021.

\bibitem[Good(1946)]{good1946normal}
I.~J. Good.
\newblock Normal recurring decimals.
\newblock \emph{Journal of the London Mathematical Society}, 1\penalty0 (3):\penalty0 167--169, 1946.

\bibitem[Heath-Brown(1992)]{heath1992density}
D.~Heath-Brown.
\newblock The density of zeros of forms for which weak approximation fails.
\newblock \emph{Mathematics of computation}, 59\penalty0 (200):\penalty0 613--623, 1992.

\bibitem[Huisman(2016)]{huisman2016newer}
S.~G. Huisman.
\newblock Newer sums of three cubes.
\newblock \emph{arXiv preprint arXiv:1604.07746}, 2016.

\bibitem[Romero et~al.(2021)]{romero2021matrixe_phi}
P.~F.~M. Romero et~al.
\newblock Matrixe\_$\phi$\^{} 9. classification of integers and primordial algebra.
\newblock \emph{MATUA. Revista de matem{\'a}ticas de la Universidad del Atl{\'a}ntico (Colombia).}, 8\penalty0 (1):\penalty0 10--45, 2021.

\end{thebibliography}

\end{document}